\documentclass[11p,leqno]{amsart}
\usepackage{amsmath}
\usepackage{amsfonts}
\usepackage{amssymb}
\usepackage{graphicx}
\usepackage{color}
\newtheorem{theorem}{Theorem}[section]\newtheorem{thm}[theorem]{Theorem}
\newtheorem*{theorem*}{Theorem}
\newtheorem{lemma}{Lemma}[section]
\newtheorem{corollary}[theorem]{Corollary}

\newtheorem{prop}{Proposition}[section]
\newtheorem{remark}[theorem]{Remark}

\def\qed{\hfill $\square$}

\def\Ric{\text{Ric}}

\def\C{\Bbb C}
\def\R{\Bbb R}
\def\CP{{\Bbb C} P}
\def\Sph{\Bbb S}
\def\Ball{\Bbb B}

\def\id{\operatorname{id}}
\def\Ric{\operatorname{Ric}}
\def\Rm{\operatorname{Rm}}
\def\I{\operatorname{I}}
\def\Scal{\operatorname{Scal}}
\def\tr{\operatorname{tr}}

\newcommand{\SO}{{\mathsf{SO}}}
\newcommand{\Oo}{{\mathsf{O}}}
\newcommand{\so}{{\mathfrak{so}}}
\def\I{\operatorname{I}}
\newcommand{\eps}{{\varepsilon}}
\newcommand{\ml}{{\langle}}
\newcommand{\mr}{{\rangle}}
\newcommand{\Spin}{{\mathsf{Spin}}}
\newcommand{\HP}{{\mathbb{H}P}}
\newcommand{\CaP}{{\mathrm{Ca}P}}

\newcommand{\vol}{\mathrm{vol}}

\numberwithin{equation}{section}

\begin{document}
\title{\bf Manifolds with $1/4$-pinched flag curvature}

\author{Lei Ni and Burkhard Wilking}



\date{}

\maketitle

\begin{abstract}
We say that a nonnegatively curved manifold $(M,g)$ has quarter pinched flag curvature
if for any two planes which intersect in a line the ratio of their sectional
 curvature is bounded above by $4$.
 We show that these manifolds have nonnegative complex sectional curvature.
 By combining with a theorem of Brendle and Schoen it follows that any positively curved
 manifold with strictly quarter pinched flag curvature must be a space form.
 This in turn generalizes a result of Andrews and Nguyen in dimension 4.
 For odd dimensional manifolds we obtain results for the case that the flag curvature
 is pinched with some constant below one quarter, one of which generalizes a
 recent work of Petersen and Tao.

\end{abstract}

\section{Introduction}

Let $(M, g)$ be a Riemannian manifold with curvature tensor $R$ and curvature operator $\Rm$ (we follow the convention $R(X, Y, Z, W)=\langle \Rm(X\wedge Y), Z\wedge W\rangle$). Assume that $(M, g)$ has nonnegative sectional curvature. Fixing a point $x\in M$, for any nonzero vector $e\in T_xM$, we  define the flag curvature in the direction $e$ by the symmetric bilinear form $R_e(X, X)=R(e, X, e, X)$. Restrict $R_e(\cdot, \cdot)$ to the subspace orthogonal to $e$, it is semi-positive definite. We call {\it  $(M, g)$ has $\lambda$-pinched flag curvature ($1>\lambda\ge 0$) if the eigenvalues of the symmetric bilinear  form $R_e(\cdot, \cdot)$, restricted to the subspace  orthogonal to $e$, are $\lambda$-pinched  for all  nonzero vector $e$}. Namely
\begin{equation}\label{eq:fcp1}
R_e(X, X)\ge \lambda (x) R_e(Y, Y)
\end{equation}
for any $X, Y$  in the subspace orthogonal to $e$,
with $|X|=|Y|$.  This condition was recently brought to the attention by the work of Andrews and Nguyen \cite{AN}.  In \cite{AN}, Andrews and Nguyen proved the following theorem.

\begin{theorem*}[Andrews-Nguyen] For any $\lambda\ge \tfrac{1}{4}$ the class of
positively curved $4$-manifolds with $\lambda$-pinched flag curvature is invariant under the Ricci flow.
Moreover any such manifold is either diffeomorphic to a spherical space form or isometric to $\CP^2$ with Fubini-Study metric (up to a  scaling).
\end{theorem*}

The $\lambda$-pinched flag curvature condition is equivalent to saying that $K(\sigma_1)\ge \lambda K(\sigma_2)$ for a pair of planes $\sigma_1$ and $\sigma_2$ such that $\sigma_1\cap \sigma_2\ne \{0\}$.
The above result generalizes a theorem of Chen \cite{Ch}, who used Hamilton's Ricci flow to classify
four manifolds with pointwise quarter pinched sectional curvature. We recall that a manifold is said to be pointwise
quarter pinched if $K(\sigma_1)\ge \tfrac{1}{4}K(\sigma_2)$ holds
 for all planes $\sigma_1$, $\sigma_2\subset T_xM$ and all $x\in M$.

  The earlier convergence results on Ricci flow are mostly for manifolds with dimension not greater than four \cite{H82,H86}.
Recently, in \cite{BW}, B\"ohm and the second author constructed a new family of cones (in the space of algebraic curvature operators) which is invariant under Ricci flow  and proved that any compact manifold with $2$-positive curvature operator is diffeomorphic to a spherical space form, by showing that the normalized Ricci flow evolves such a metric to one with constant curvature.

Soon afterwards, Brendle and Schoen \cite{BS} proved that any manifold with strictly quarter pinched sectional curvature is a spherical space form. The key novel step in the proof is to show that nonnegative isotropic curvature is invariant under the Ricci flow,  which was also independently proved by Nguyen \cite{Ng}. Another important step was to show that for any pointwise $1/4$-pinched manifold $M$, the manifold $(M, g)\times \R^2$ has nonnegative isotropic curvature. The convergence of the normalized Ricci flow for metrics with strictly $1/4$-pinched sectional curvature then follows from a convergence result in \cite{BW} (cf. Theorem 3.2 in Section 3). Later on it was pointed out in \cite{NW} that $(M, g)\times \R^2$ has nonnegative isotropic curvature if and only if $(M, g)$ has nonnegative complex sectional curvature. (The readers should consult \cite{BW}, \cite{BS}, and  Section 2 for notations involved.) The convergence result of \cite{BS} can then be restated as

\begin{theorem*}[Brendle-Schoen]
A metric $g$ of positive complex sectional curvature on a compact manifold evolves under the normalized Ricci flow to a constant curvature limit metric.
\end{theorem*}

The part that strictly pointwise $1/4$-pinched manifolds have positive complex sectional curvature was essentially  known to be true by an argument of \cite{YZ} and \cite{He} for the negatively $1/4$-pinched manifolds. The main result of this paper is to generalize the latter result by relaxing the assumption of pointwise $1/4$-pinched sectional curvature to the assumption of  $1/4$-pinched flag curvature.

\begin{theorem} \label{main1} Let $(M^n, g)$ be a nonnegatively curved  Riemannian manifold.
 If $(M,g)$ has quarter pinched flag curvature, then
$(M, g)$ has nonnegative complex sectional curvature.
\end{theorem}

This is proved in Section 2.
Combining the last two results and the rigidity theorem of  \cite{BS2} on manifolds with weakly $1/4$-pinched sectional curvature, we obtain the following corollary as a generalization.
\begin{corollary}\label{weak-fqp} Let $(M^n, g)$  be a compact nonnegatively curved Riemannian mani\-fold with $1/4$-pinched flag curvature and the scalar curvature $\Scal(x)> 0$ for some $x\in M$.  Then $(M, g)$ is diffeomorphic to a spherical space form or isometric to finite quotient of a rank one symmetric space.
\end{corollary}
This is done in Section 3.
We should point out that if an algebraic curvature operator has $\lambda$-pinched  flag curvature, then its sectional curvature is $\lambda^2$-pinched.
We will see in Section 4 that this inequality is indeed sharp in dimensions above 3.

A result similar to Theorem \ref{main1} yields the following generalization of a result by Yau and Zheng \cite{YZ}, and Hern\'andez \cite{He}, by relaxing the condition on $1/4$-pinched sectional curvature to $1/4$-pinched flag curvature.

\begin{theorem}\label{main2}
If $(M^m, g)$ ($m\ge 2$) is a compact  complex K\"ahler manifold which also admits a Riemannian metric $h$ with negative $1/4$-pinched flag curvature. Then $(M, h)$ must be holomorphically isometric to a compact quotient of the unit ball $\Ball^m\subset \C^m$.
\end{theorem}

This is also done in Section 4. In Section 5 we consider nonnegatively curved manifolds of odd dimension $n=2m+1$ with $\lambda$-pinched flag curvature for $\lambda<\frac{1}{4}$. A generalization of an earlier result of Berger \cite{Berger} on the vanishing of the second Betti number is obtained for $\lambda\ge\frac{n-3}{4n-9}$ (cf. Theorem \ref{odd-berger}), provided that the scalar curvature $\Scal(x)>0$ for some $x\in M$.  An interesting  corollary  is for $5$ dimensional manifolds.

\begin{corollary}\label{berger-5d} If $(M, g)$ is a closed positively curved $5$-manifold with  $\frac{2}{11}$-pinched flag curvature, then $M$ must be a rational-homological sphere.
\end{corollary}

An  open question is whether or not a $5$-manifold as above is diffeomorphic to a spherical space form.
For odd dimensional manifolds, we also proved a result generalizing an earlier one of Petersen and Tao \cite{peter-tao} for nonnegatively curved manifolds (even orbifords) with flag pinching constant below one quarter and an arbitrarily small  global curvature pinching condition.

\begin{thm} For any  odd number $n=2m+1$ and $C>0$,  there is an $\eps>0$
such that the following holds.
Let $(M^n, g)$ be a nonnegatively curved Riemannian  orbifold  of dimension $n$ with $\tfrac{1-\eps}{4}$ pinched-flag curvature and scalar curvature satisfying $1\le \Scal \le C$. Then $M$ admits a metric of constant curvature.
\end{thm}

This was done in Section 6. The proof makes use of the above mentioned generalization (of an earlier  result of Berger) on the vanishing of the second Betti number as well as the injectivity radius estimates of Petrunin-Tuschmann and Fang-Rong \cite{PTu, FR}).  We also explain
that  in even dimensions for $\eps>0$ there are infinitely many distinct simply connected
 orbifolds (namely certain weighted projective space) with
curvature $\tfrac{1-\eps}{4}<K<1$.  Nevertheless
we are able to classify orbifolds
with nearly quarter flag pinching in even dimensions as well, see Theorem \ref{orbi2}.
This might be interesting in other context as well since its proof provides, in very special circumstances,
 a method to deal with a sequence of collapsing solutions of the Ricci flow.

\section{Pinched flag curvature and  the complex sectional curvature}

Fix a point $x\in M$. Let $T^{\C}_xM\doteqdot T_xM\otimes \C$ be the complexified tangent bundle. The curvature operator,  viewed as a symmetric tensor (transformation) of $ \wedge^2 T_xM$, can be extended to $\wedge^2 T^{\C}_xM$. We also extend the inner product $\langle \cdot, \cdot \rangle$ bilinearly to $T^{\C}_xM$.
We defined the complex sectional curvature of a pair of vectors $U, V$ by $R(U, V, \bar{U}, \bar{V})\doteqdot \langle \Rm (U\wedge V), \overline {U \wedge V}\rangle$. This notion of curvature appeared previously in geometry mainly in the study of harmonic maps and its applications to the rigidity theorems (cf. \cite{Sa}, \cite{Siu}). The interested reader may also  consult \cite{Gr} for its relations with various other curvature notions. In particular, we say that $(M, g)$ has nonnegative isotropic curvature if $R(U, V, \bar{U}, \bar{V})\ge 0$ for any $U$ and $V$ which spans an isotropic (with respect to $\langle \cdot, \cdot\rangle$)  plane. If $U=X+\sqrt{-1}Y, \quad V=Z+\sqrt{-1}W$, via the first Bianchi identity, we have that
\begin{eqnarray}
R(U, V, \bar{U}, \bar{V})&=& R(X, Z, X, Z)+R(Y, W, Y, W)+R(X, W, X, W)+R(Y, Z, Y, Z)\nonumber\\
&\quad& -2R(X, Y, Z, W). \label{eq:cinr}
\end{eqnarray}
The aim of this section is to show that $\lambda$-pinched flag curvature implies the positivity of the complex sectional curvature for $\lambda>\frac{1}{4}$. The proof relies on an improved version of Berger's lemma \cite{Berger} and a modification of the argument of \cite{YZ}. (See also \cite{He, NW}.)

First we start with a simple known lemma. We include its proof for the sake of the completeness.
\begin{lemma}
\begin{eqnarray}
6R(X, Y)Z&=& -R(Y, Z+X) (Z+X) +R(Y, Z-X)(Z-X)
\nonumber\\ &&+R(X, Z+Y) (Z+Y)-R(X, Z-Y)(Z-Y).\label{eq:polar1}
\end{eqnarray}
\end{lemma}
\noindent {\sl Proof.} First Bianchi yields:
$$
-R(X, Y)Z=R(Z, X)Y+R(Y, Z)X.
$$
While
\begin{eqnarray*}
R(Y, Z)X &=& \frac{1}{2}\left(R(Y, Z+X)(Z+X)-R(Y, Z-X)(Z-X)\right)-R(Y,X)Z;\\
R(Z,X)Y &=& -R(X, Z)Y\\
&=& -\frac{1}{2}\left( R(X, Z+Y)(Z+Y)-R(X, Z-Y)(Z-Y)\right)+R(X, Y)Z.
\end{eqnarray*}
Combining them we have the claimed equation.
\qed

\begin{corollary}\label{polar2} Let $k(A, B)\doteqdot R(A, B, A, B)$.
\begin{eqnarray*}
-R(X, Y, Z, W)= \frac{1}{12}\left(\begin{array}{l}-k(Y+W, Z+X) +k(Y, Z+X)+k(W, Z+X)\\
-k(X+W, Z-Y)+k(X, Z-Y)+k(W, Z-Y)\\
 -k(Y-W, Z-X)+k(Y, Z-X)+k(W, Z-X)\\
-k(X-W, Z+Y)+k(X, Z+Y)+k(W, Z+Y)\end{array}\right).
\end{eqnarray*}
\end{corollary}
\noindent {\sl Proof.}
By the lemma we have that
\begin{eqnarray*}
-R(X, Y, Z, W)&=&\frac{1}{6}\left( R(Y, Z+X, Z+X, W) -R(Y, Z-X, Z-X, W)\right.
\\&\quad&\left.-R(X, Z+Y, Z+Y, W)+R(X, Z-Y, Z-Y, W)\right).
\end{eqnarray*}
The result follows from the fact that $R( \cdot, Z+X, Z+X, \cdot)$, $R( \cdot, Z-X, Z-X, \cdot)$, etc are symmetric bilinear forms and the standard polarization formula.
\qed

The next lemma is the key to estimate the complex sectional curvature in terms of the pinching constant of  the flag curvatures.

\begin{lemma} \label{help1}Let $e\in T_xM$ be a nonzero vector. Assume that the sectional curvature is nonnegative at $x$. Let $b(\cdot, \cdot) \doteqdot R_e(\cdot, \cdot)$, and let $\mathcal{N}$ be the subspace orthogonal to $e$. Assume that the flag curvature of $e$ is $\lambda$-pinched with $\lambda>0$. Then for any $Y$,  $W$  $\in \mathcal{N}$ such that $\langle Y, W\rangle =0$,
\begin{equation}\label{eq:key-est1}
\frac{2\lambda}{1+\lambda}(b(Y, Y)+b(W, W))\le b(Y+W, Y+W)\le \frac{2}{1+\lambda}(b(Y, Y)+b(W, W)).
\end{equation}
 If equality holds and $b\neq 0$, then  $|W|=|Y|$.
\end{lemma}
\noindent {\sl Proof.} Let $a=\frac{1}{\lambda}$. We may restrict $b(\cdot, \cdot)$ to the $2$-plane $\mathcal{P}_2$ spanned by $Y, W$. Assume that
$b(e_1, e_1)=\max_{Z\in \mathcal{P}_2, |Z|=1}  b( Z, Z)=\hat{\Lambda}$ and $b(e_2, e_2)=\min_{Z\in \mathcal{P}_2, |Z|=1}  b( Z, Z)=\hat{\lambda}$. Without loss of the generality we can assume that $\hat{\lambda}>0$. Let $\hat{a}=\frac{\hat{\Lambda}}{\hat{\lambda}}$. Since one can multiply both sides of (\ref{eq:key-est1}) by a constant, without loss of the generality we can assume that $|Y|=1$ and $|W|\le 1$. Write $Y=\cos \theta e_1 + \sin \theta e_2$ and $W=\eta(-\sin \theta e_1 +\cos \theta e_2)$ with $|\eta|<1$.
For the upper bound we just compute $b(Y+W, Y+W)$ and $b(Y, Y)+b(W, W)$ and compare them. Note also $\hat{a}\le a$. The computation shows that the claimed upper bound in (\ref{eq:key-est1}) is, after dividing $\hat{\lambda}$ on both side,   equivalent to
\begin{eqnarray*}
(\hat{a}-1) \cos^2 \theta &+&(\hat{a}-1) \eta^2 \sin^2 \theta +1+\eta^2 -2(\hat{a}-1)\eta \sin \theta \cos \theta \\
&\le& \frac{2a}{a+1}\left( 1+\eta^2+(\hat{a}-1)\cos^2 \theta +\eta^2 (\hat{a}-1) \sin^2\theta\right).
\end{eqnarray*}
This can be further reduced to
$$
-2(\hat{a}-1)\eta \sin \theta \cos\theta \le \frac{a-1}{a+1}(\eta^2+1) +\frac{a-1}{a+1} (\hat{a}-1)\cos^2\theta +\frac{a-1}{a+1}(\hat{a}-1)\eta^2 \sin ^2 \theta.
$$
Using that $a\ge \hat{a} \ge 1$, it suffices to show that
$$
-2(a+1)\eta \sin\theta \cos \theta \le (1+\eta^2) +(a-1)\cos^2\theta +\eta^2 (a-1)\sin^2 \theta.
$$
This last inequality is elementary, noting that $a\ge1$. The proof of the lower bound is similar. Clearly the equality does imply that $|\eta|=1$, namely $|Y|=|W|$.
\qed

The consequence of the above lemma is a generalization of Berger's well-known lemma, which we state below.

\begin{corollary} \label{Berger} Assume that $(M, g)$ has $\lambda$-pinched flag curvature with dimension $n\ge 4$. Assume that the sectional curvature is nonnegative at $x$ and $X, Y, Z, W \in T_x M $ are four vectors mutually orthogonal.
  Then
\begin{eqnarray*}
 6\frac{1+\lambda}{1-\lambda}|R(X, Y, Z, W)|&\le &k(Y, Z)+k(X, Z) +k(X, W)+k(Y, W)
+2k(X, Y)+2k(Z, W).
\end{eqnarray*}
 If equality holds and $\Rm (x)\neq 0$, then  vectors $X, Y, Z, W$  have the same norm.
\end{corollary}
\noindent {\sl Proof.} As before let $a=\frac{1}{\lambda}$. The result follows from Corollary \ref{polar2} and the estimate (\ref{eq:key-est1}). In fact, applying (\ref{eq:key-est1}) one can estimate
\begin{eqnarray*}
&\,&R(Y+W, Z+X, Z+X, Y+W)-R(Y, Z+X, Z+X, Y)-R(W, Z+X, Z+X, W)\\
&=&  R(Y, Z+X, Y, Z+X)+R(W, Z+X, W, Z+X)-R(Y+W, Z+X, Y+W, Z+X)\\
&\quad& \quad \quad \ge \left(1-\frac{2a}{a+1}\right)\left(R(Y, Z+X, Y, Z+X)+R(W, Z+X, W, Z+X)\right).
\end{eqnarray*}
Applying the similar estimate to the other three groups in the expression of $R(X, Y, Z, W)$ from Corollary \ref{polar2}, the lower bound of $-R(X, Y, Z, W)$
follows by expanding $R(Y, Z+X, Y, Z+X)+R(W, Z+X, W, Z+X)$, etc  and regrouping them. The upper bound is similar. We leave the detailed checking to the readers.
\qed

\begin{remark} The Berger's lemma concludes a similar estimate as that of  Corollary \ref{Berger} under the stronger assumption that the sectional curvature is $\lambda$-pinched.
\end{remark}

Next we shall estimate the complex sectional curvature via the pinching constant $\lambda$.

\begin{lemma}\label{lemma-lb-cs} Assume that $(M, g)$ (with dimension $n\ge 4$) has nonnegative sectional curvature  and satisfies (\ref{eq:fcp1}) at $x\in M$. Then for any $X, Y, Z, W \in T_x M $ which are mutually orthogonal with
$|X|=|Z|=1$, $1\ge |Y|\ge |W|\ge 0$,  we have the following estimates
\begin{eqnarray*}
R(U, V, \bar{U}, \bar{V})
&\ge& \frac{2(4\lambda-1)}{3}\left(k(X, \hat{Y})+k(Z, \hat{W})\right)
\end{eqnarray*}
when $\lambda\le \frac{1}{4}$, and
$$
R(U, V, \bar{U}, \bar{V})
\ge \frac{2(4\lambda-1)}{3}\left(k(X, \hat{Y})\right)
$$
when $\lambda>\frac{1}{4}$.
Here $U=X+\sqrt{-1}Y$, $V=Z+\sqrt{-1}W$,  $\hat{Y}$ and $\hat{W}$ are the unit vectors in the direction of $Y$ and $W$. The equality in the case $\lambda\le \frac{1}{4}$  can only happen when $|X|=|Y|=|W|=|Z|$ and $\lambda k(X, Y)=\lambda k(Z, W)=k(X, Z)=k(X, W)=k(Y, Z)=k(Y, W)$, provided that the sectional (or scalar) curvature  is positive.
\end{lemma}
\noindent {\sl Proof.} Recall that $R(U, V, \bar{U}, \bar{V})=k(X, Z)+k(X, W)+k(Y, Z)+k(Y, W)-2R(X, Y, Z, W)$. By Corollary \ref{Berger} we have that
\begin{eqnarray*}
R(U, V, \bar{U}, \bar{V})
&\ge& \frac{2(a+2)}{3(a+1)}\left(k(X,Z)+k(X, W)+k(Y, Z)+k(Y, W)\right)\\
&\quad&-\frac{2(a-1)}{3(a+1)}\left(k(X, Y)+k(Z, W)\right).
\end{eqnarray*}

Now write $t=|Y|^2, \xi =|W|^2$ and $Y= t^{1/2} \hat{Y}, \quad W= \xi^{1/2} \hat{W}$. Consider
\begin{eqnarray*}
Q(t, \xi)&=&(a+2)\left(k(X, Z)+ k(X, \hat{W}) \xi +k(\hat{Y}, Z) t+k(\hat{Y}, \hat{W}) t\xi\right)\\
&\quad&-(a-1)\left(k(X, \hat{Y})t+k(Z, \hat{W}) \xi\right).
\end{eqnarray*}
as a function defined for $0\le \xi\le t \le 1$. Its minimum can only possibly achieve  when
$\xi=t=0$, $\xi=0, t=1$ or $\xi=t=1$. For the last case we have that
$$
Q\ge \frac{(4-a)(a+1)}{a}\left(k(X, \hat{Y})+k(Z, \hat{W})\right)
$$
which implies the result. The other two cases are similar (easier).
\qed

\noindent {\sl Proof of Theorem \ref{main1}.} Let $\Rm\in S^2_B(\so(n))$ be an
algebraic curvature operator with quarter pinched flag curvature.
We want to check that the curvature of $\Rm$ evaluated at a complex plane
 $\sigma\subset \mathbb{C}^n$ is nonnegative.
Next we choose  a good basis of $\sigma$ as follows.
We endow $\mathbb{C}^n$ with the usual scalar product
$\ml \ml x,y\mr\mr= \bar{y}^{tr}x$.
Let $U$ be a unit vector in $\sigma$
which maximizes among all unit vectors the real part of the scalar product
$\langle\ml U,\bar{U}\mr\mr$.
Let $V$ be a unit vector in $\sigma$ perpendicular to
$U$ such that $\ml \ml V,\bar{V}\mr\mr$ is a nonnegative real number.
By the choice of $U$ we have $\ml\ml U,V\mr\mr=\ml \ml\bar{U},V\mr\mr=0$.
Moreover, $\ml \ml U, \bar{U}\mr\mr$  and $\ml\ml V,\bar V\mr\mr$ are  real.
Let $X,Y,Z,W\in \mathbb{R}^n $ be such that
$U=X+\sqrt{-1}Y$, $V=Z+\sqrt{-1}W$.
It is straightforward to check that $X,Y,Z,W$ are pairwise orthogonal.
Furthermore $|X|\ge|Y|$ and $|Z|\ge |W|$.

Finally we rescale $U$ and $V$ by a real number such that
$|X|=|Z|=1$. The nonnegativity now follows from the previous Lemma.

\qed

The above argument  also proves that
a curvature operator $\Rm\neq 0$ with nonnegative sectional curvature,
$\lambda$-pinched flag curvature has positive complex sectional curvature if $\lambda>\tfrac{1}{4}$. Tracing the proof it is easy to see the following corollary.

\begin{corollary}\label{equality} Let $(M, g)$ be as in Theorem \ref{main1}. Then if $R(U, V, \bar{U}, \bar{V})=0$ for a pair of vectors  $U, V$ spanning a complex plane $\sigma$. Then $\sigma$ is  isotopic. Namely for any $T\in \sigma$, $\langle T, T\rangle =0$.
\end{corollary}

\section{Consequences}
Using the work of \cite{BW} and \cite{BS}, Theorem \ref{main1}
 has a more general  differential sphere theorem as its direct consequence. Recall that it was proved in \cite{BS} that

\begin{theorem}[Brendle-Schoen]
\label{thm:PCSC}
The Ricci flow on a compact manifold preserves the cones
consisting of: (i) the curvature operators with nonnegative complex
sectional curvature and (ii) the curvature operators with positive complex
sectional curvature.
\end{theorem}

Originally in \cite{BS}, the authors proved that the condition that $M \times \R^2$ has nonnegative  isotropic curvature is preserved under Ricci flow. However it turns out (cf. \cite{NW}) that the condition $M\times \R^2$ has nonnegative isotropic curvature is equivalent to $M$ has nonnegative complex sectional curvature. There  exists a more direct proof of the above Brendle-Schoen's result in \cite{NW} using  complex numbers.

 Using the above result, Brendle and Schoen further proved that any compact manifold with strictly $1/4$-pinched sectional curvature must be a diffeomorphic sphere, via the following general invariant cone construction of \cite{BW} by B\"ohm and the second author, and their consequence on the convergence of the Ricci flow. (Cf. Theorem 6.2 of \cite{W}.)

\begin{theorem}[B\"ohm-Wilking]\label{thm:RF-con} Let $C$ be an $\Oo(n)$ invariant cone in the vector space of curvature operators $S^2_B(\mathfrak{so}(n))$ with the following properties

(i) $C$ is invariant under the ODE $\frac{d \Rm}{dt} =\Rm^2+\Rm^{\#}$.

(ii) $C$ contains the cone of the nonnegative curvature operator or slightly weaker all nonnegative curvature operators of rank $1$.

(iii) $C$ is contained in the cone of curvature operators with nonnegative sectional curvature.

Then for any compact manifold $(M, g)$ whose curvature operator is contained in the interior of $C$ at every point $p\in M$, the normalized Ricci flow evolves $g$ to a limit metric of constant sectional curvature.
\end{theorem}
\begin{proof} It follows from the proof of Lemma 3.4 and 3.5 of \cite{BW} verbatim.
\end{proof}
Combining Theorem \ref{thm:PCSC} and Theorem \ref{thm:RF-con} we have that if $(M, g)$ is a Riemannian manifold with positive sectional curvature such that (\ref{eq:fcp1}) holds for some $\lambda> \frac{1}{4}$,  then $(M, g)$ is diffeomorphic to a spherical space form. In fact one can arrive the same conclusion assuming only that $(M,g)$ is
a nonnegatively curved manifold with quarter pinched flag curvature and with some point $x$ at which the
curvature operator $\Rm(x)\neq 0$ has strictly quarter pinched flag curvature.

In \cite{BS2}, the manifolds with weakly $1/4$-pinched sectional curvature was analyzed. Since the key step is a result on the rigidity of manifolds with nonnegative complex sectional curvature (Theorem 2 of \cite{BS2}), which in turn relies on a general strong maximum principle result,  one can modify this part  to obtain the generalization stated in Corollary  \ref{weak-fqp}.

\noindent {\sl Proof of Corollary \ref{weak-fqp}.} By \cite{BS2}, particularly, Proposition 11 therein on the rigidity of manifolds with nonnegative complex sectional curvature, one can reduce the proof of Corollary \ref{weak-fqp} to show that {\it if $M$ is a K\"ahler manifold with $\frac{1}{4}$-pinched flag curvature, then the universal cover of $M$ is isometric to $\CP^m$ ($n=2m$)}. Now consider the complex sectional curvature $R(U, V, \bar{U}, \bar{V})$ for any $U, V \in T^{1,0}M$ with $|U|=|V|$ and $\langle U, \bar{V}\rangle =0$. By the K\"ahlerity it must be zero. Write as before that $U=X+\sqrt{-1}Y$ and $V=Z+\sqrt{-1}W$. By tracing the equality case in Section 2, we deduce that $Y=-JX$, $W=-JZ$ and $k(X, JX)=k(Z, JZ)=4k(X,Z)=4k(X, W)=4k(Y, Z)=4k(Y, W)$. This proves the result for the case $m\ge 3$ since for any $U_1, U_2 \in T^{1,0}M$ one can find $U_3$ such that $\langle U_3, \bar{U}_i\rangle =0$ for $i=1, 2$. Hence $M$ has constant holomorphic sectional curvature at any given point $p\in M$. By the K\"ahlerian analogue of Schur's lemma we can conclude that $M$ is of constant holomorphic sectional curvature.  For $m=2$, the result has been proved in \cite{AN}. (See also the next section for a different argument on this part.)
\qed

\section{An example and the dual case}

In this section, first for $n\ge 4$,  we present an example of algebraic curvature operator which has $1/4$-pinched sectional curvature, but only $1/2$-pinched flag curvature.

Define an algebraic curvature operator on $\R^4$ by
\begin{eqnarray*}\Rm &=&4(e_1\wedge e_2)\otimes (e_1\wedge e_2)+2(e_1\wedge e_3)\otimes (e_1\wedge e_3)\\
&+&2(e_1\wedge e_4)\otimes (e_1\wedge e_4)+2(e_2\wedge e_3)\otimes (e_2\wedge e_3)+2(e_2\wedge e_4)\otimes (e_2\wedge e_4)\\&+&(e_3\wedge e_4)\otimes (e_3\wedge e_4).
\end{eqnarray*}
Here $\{e_i\}$ is a orthonormal frame. First one can easily check that $\Rm$ is an algebraic curvature operator since it is a symmetric tensor of $\wedge^2(\R^4)$ and satisfies the first Bianchi identity.
It is easy to see that the sectional curvature of $\Rm$ is $1/4$-pinched. On the other hand, for the flag curvature it is only $1/2$-pinched. It is  clear that the eigenvalues of $R_e(\cdot, \cdot)$ are $1/2$-pinched, when $e=e_i$. To check for generic vector it is sufficient to  check for $e=\cos\theta e_1+\sin \theta e_3$ for any $\theta\in \R$. Direct computation verifies that the three eigenvectors of $R_e(\cdot, \cdot)$ are
$e_2, e_4 $ and $-\sin \theta e_1+\cos \theta e_3$ with eigenvalues $2\cos^2\theta+2$, $1+\cos^2\theta$ and $2$ correspondingly. They have pinching constant $1/2$ again.
The example shows that if an algebraic curvature tensor $R$ has $\lambda$-pinched flag curvature, its sectional curvature  has pinching constant at the best $\lambda^2$. In particular, if $\Rm$ has $\frac{1}{4}$-pinched flag curvature, one at best can expect its sectional curvature is $\frac{1}{16}$-pinched. This example can easily be generalized to any dimensions by replacing $e_3$ and $e_4$ by $e_i$ for $n\ge i\ge 3$. One may be able to endow the sphere of dimension above 3 with a metric such that its curvature operator approximates, , after scaling,  the above algebraic curvature operator pointwise with any given precision.

Finally we observe that the discussions in Section 2 can be adapted to the manifolds with nonpositive sectional curvature and $1/4$-pinched flag curvature. Note that in the case that $(M,h)$ has nonpositive sectional curvature we say that  it has $\lambda$-pinched flag curvature if and only if
\begin{equation}\label{eq:fcp2}
\lambda R_e(X, X)\ge R_e(Y, Y)
\end{equation}
holds for any $e\ne 0$, and any $X, Y$  in the orthogonal complement of $e$ with $|X|=|Y|=1$.
Clearly this is equivalent to saying that the curvature tensor $R'=-R$ has nonnegative sectional
curvature with $\lambda$-pinched flag curvature.
The following is the analogue of Theorem \ref{main1}

\begin{corollary}\label{negative-flag} Assume that  $(M^n, g)$ is a Riemannian manifold ($n\ge 4$) having nonpositive sectional curvature. If (\ref{eq:fcp2}) holds for some $\lambda\ge \frac{1}{4}$, then
$(M, g)$ has nonpositive complex sectional curvature. Moreover if the sectional curvature is negative at $x$ and (\ref{eq:fcp2}) holds at $x$ with $\lambda>\frac{1}{4}$, for all $e(x)\ne0$, then there exists $\epsilon(x)>0$, depending on $\lambda(x)$ and $\Scal(x)$ such that at $x$, $(\Rm +\epsilon \I)$ has nonpositive complex sectional curvature, where $\I$ is the identity of $S^2(\wedge^2(\R))$.
In particular $(M, g)$ has negative complex sectional curvature at $x$.
\end{corollary}

Now Theorem \ref{main2} follows from the existence result of Eells and Sampson on the harmonic maps and the following general rigidity theorem.

\begin{theorem} \label{sampson-1}Let $X^m$ ($m\ge 2$) be a compact K\"ahler manifold of complex dimension $m$, $M^n$ be a Riemannian manifold of even  dimension (real) $n=2m$ with non-positive sectional curvature and satisfying (\ref{eq:fcp2}) with $\lambda=\frac{1}{4}$. Suppose that there is continuous map
 $f_0: X\to M$ of nonzero degree. Then $M$ is a locally symmetric space with universal cover $\Ball^m$ ($n=2m$).
\end{theorem}
\noindent {\sl Proof.} Here we follow the line of argument in \cite{YZ}.   First by the existence theorem of Eells and Sampson $f_0$ is homotopic to a harmonic map $f$ which is onto. Let $\mathcal{Z}$ be the critical value of $f$. First restrict the discussion on $M\setminus \mathcal{Z}$. Endow $M \setminus \mathcal{Z}$ with an almost complex structure (still denoted by $J$) by pushing forward the complex $J$ on $X$ via $f_*$. Let $(z_1, \cdot\cdot \cdot, z_m)$ be a local complex coordinate of $X$. Let   $U_i=f_*(\frac{\partial}{\partial z_i})$. Since by Corollary \ref{negative-flag} $M$ has non-positive complex sectional curvature, a result of Sampson, Theorem 1 of  \cite{Sa}, asserts that
$$
\nabla_{U_i} \overline{U}_j =0, \quad \quad R(U_i, U_j, \overline{U}_i, \overline{U}_j)=0.
$$
This particularly implies that $R(U, V, \bar{U}, \bar{V})=0$ for any $U, V \in T^{1,0} M$ (here the decomposition is with respect to the push-forward complex structure). The dual version of Corollary \ref{equality} asserts that $\sigma=\mbox{Span}\{U, V\}$ is isotropic, if it is a complex plane.  This in turn implies that $\langle U, U\rangle=0$, hence the metric $h$, with respect to the push-forward of the complex structure, is Hermitian. Using $\nabla_{U_i}\overline{U}_j=0$,  one can check, via a direct computation,  that the K\"ahler form of  $(M\setminus \mathcal{Z}, h, J)$  is closed (this is the observation of \cite{YZ}, the proof of Theorem 3). The K\"ahlerity conclusion is first made on $M \setminus \mathcal{Z}$  and then extends to $M$ since $M \setminus \mathcal{Z}$ is open and dense. Once we know that $(M, h, J)$ is a K\"ahler manifold with $1/4$-pinched flag curvature, when $m\ge 3$, the claimed result follows similarly as in the proof of Corollary \ref{weak-fqp}.  More precisely, we can show,  in exactly the same manner as in Corollary \ref{weak-fqp} that, at every point $x\in M$, its holomorphic sectional curvature is  independent of the choice of the complex lines. Hence the result follows by the holomorphic  Schur's lemma. When $m=2$, the last conclusion can be seen by observing that at every point $x\in M$, the curvature operator $\Rm$ is the multiple of curvature operator of  $\Ball^2$. The reason is that the condition of $1/4$-flag curvature pinching is a convex $U(2)$ invariant condition including the curvature operator of $\Ball^2$. On the other hand, one can check easily that there exist some small perturbations towards the holomorphic Weyl part, as well as the  traceless Ricci part, which  violates  the condition of the  flag curvature being $1/4$-pinched. This shows that the curvature operator of $\Ball^2$ is the only one curvature operator of K\"ahler manifolds satisfying the flag curvature $1/4$-pinching condition.
\qed

\section{Odd dimensions--a vanishing theorem}

As another application of the estimates from Section 2 we prove a generalization of an earlier result of Berger (Theorem 2 of \cite{Berger}) on odd dimensional manifolds. Again, we relax the assumption on the pointwise $\lambda$-pinching of the sectional curvature  to the $\lambda$-pinched flag curvature. The difference on the assumption does cause subtlety which needs to be dealt with more carefully.

\begin{theorem}\label{odd-berger}  Let $(M, g)$ be a nonnegatively curved compact Riemannian manifold of odd dimension $n$. Assume that it has  $\lambda$-pinched flag curvature with $\lambda\ge \frac{n-3}{4n-9}$ and its   scalar curvature satisfying $\Scal(x)>0$ for some point $x\in M$. Then the second Betti number of $M$ vanishes.
\end{theorem}
\noindent {\sl Proof.} Let $n=2m+1$. As in \cite{Berger} the proof is via the Bochner-Weitzenb\"ock formula:
 $$
\Delta_d=-\tr \nabla ^2 +\sum_{i, j} \omega_i \wedge i(E_j)\mathcal{R}_{E_i E_j}
$$
 on harmonic $2$-forms. Here $E_i$ is a orthonormal frame, $\mathcal{R}$ is the curvature of the induced covariant differentiation on two forms, $i(\cdot)$ is the contraction operator. We first prove the result for the case that  $\lambda> \frac{2(m-1)}{8m-5}$.
 By Hodge theory, it suffices to show that any harmonic $2$-form is trivial. With some  calculations one can show that for any two form $\Omega$
 \begin{equation}\label{help-sharp}
 \sum_{i, j} \omega_i \wedge i(E_j)\mathcal{R}_{E_i E_j} \Omega =(\Ric\wedge \id -\Rm)(\Omega).
 \end{equation}
 (Here $A\wedge B(X\wedge Y)\doteqdot \frac{1}{2}\left(A(X)\wedge B(Y)+B(X)\wedge A(Y)\right)$.)
The key is to show that if  $\lambda> \frac{2(m-1)}{8m-5}$, then $\Ric\wedge \id -\Rm >0$ as a symmetric tensor on $\wedge^2(\R^n)$ (which can be identified with $\mathfrak{so}(n)$). Also observe an interesting fact that by Lemma 2.1 of \cite{BW},  $\Ric\wedge \id -\Rm=\Rm \# \I$, even though we do not make use of it in our proof. Here $\id: \R^n \to \R^n $, $\I: S_B^2( \mathfrak{so}(n))\to S_B^2( \mathfrak{so}(n))$ are two identity maps.

Given  $\Omega\in \wedge^2(\R^n)$, by  linear algebra one can find orthonormal frame $X_1, X_2, \cdot \cdot \cdot, X_{2m+1}$ such that
$$
\Omega=\sum_{i=1}^m \alpha_i X_{2i-1}\wedge X_{2i}.
$$
Computation shows that
\begin{eqnarray*}
2\langle (\Ric\wedge \id -\Rm)(\Omega), \Omega\rangle &=& \sum_{1\le j\le m, 1\le s\le n} \alpha_j^2\left( k(2j-1, s)+k(2j, s)\right)\\
&\quad&-2\sum_{1\le i,j\le m}\alpha_i \alpha_j R(X_{2i-1}, X_{2i}, X_{2j-1}, X_{2j}).
\end{eqnarray*}
Here $k(2j, s)\doteqdot k(X_{2j}, X_s)$. Writing the last term above as
$$2\sum_{1\le j\le m} \alpha_j^2 k(2j-1, 2j)+2\sum_{1\le i\ne j\le m}\alpha_i \alpha_j R(X_{2i-1}, X_{2i}, X_{2j-1}, X_{2j})$$
and combining the similar terms we arrive at
 \begin{eqnarray*}
&\quad&2\langle (\Ric\wedge \id -\Rm)(\Omega), \Omega\rangle\\
 &=& \frac{1}{2}\sum_{i\ne j} (\alpha_i^2+\alpha_j^2) \left(k(2j-1, 2i)+k(2j-1, 2i-1)+k(2j, 2i)+k(2j, 2i-1)\right)\\
 &\,&+ \sum_{1\le j\le m} \alpha_j^2 \left(k(2j-1, n)+k(2j, n)\right)\\
 &\,&-2\sum_{1\le i\ne j\le m}\alpha_i \alpha_j R(X_{2i-1}, X_{2i}, X_{2j-1}, X_{2j}).
\end{eqnarray*}
Here $\sum_{i\ne j}$ means $\sum_{1\le i \ne j \le m}$.
Now we use Corollary \ref{Berger} to estimate the last term. Let
$$
\beta=\frac{3(1+\lambda)\lambda}{2(1-\lambda)(m-1)}
$$
and write
\begin{eqnarray*}
2|\alpha_i| |\alpha_j| R(X_{2i-1}, X_{2i}, X_{2j-1}, X_{2j})&=&2(1-\beta)|\alpha_i| |\alpha_j|R(X_{2i-1}, X_{2i}, X_{2j-1}, X_{2j})\\
&+& 2\beta\, R(|\alpha_i|^{\frac{1}{2}}X_{2i-1}, |\alpha_i|^{\frac{1}{2}}X_{2i}, |\alpha_j|^{\frac{1}{2}}X_{2j-1}, |\alpha_j|^{\frac{1}{2}}X_{2j}).
\end{eqnarray*}
Now via Corollary \ref{Berger} the right hand side above can be estimated by
\begin{eqnarray*}
&\,&\frac{1-\lambda}{3(1+\lambda)} |\alpha_i||\alpha_j|\left(k(2i, 2j-1)+k(2i, 2j)+k(2i-1, 2j-1)+k(2i-1, 2j)\right)\\
&+&\frac{2(1-\lambda)}{3(1+\lambda)}(1-\beta)|\alpha_i||\alpha_j|\left(k(2i, 2i-1)+k(2j, 2j-1)\right)\\
&+&\frac{2(1-\lambda)}{3(1+\lambda)} \beta \left( k(2i, 2i-1)\alpha_i^2+k(2j, 2j-1)\alpha_j^2\right).
\end{eqnarray*}
Observe that
\begin{eqnarray*}
 \sum_{1\le j\le m} \alpha_j^2 \left(k(2j-1, n)+k(2j, n)\right)-\sum_{i\ne j}\frac{2(1-\lambda)}{3(1+\lambda)} \beta \left( k(2i, 2i-1)\alpha_i^2+k(2j, 2j-1)\alpha_j^2\right)\ge 0
\end{eqnarray*}
by the choice of $\beta$ and the $\lambda$-pinching  of the flag curvature. Putting together we have that
\begin{eqnarray*}
&\,&2\langle (\Ric\wedge \id -\Rm)(\Omega), \Omega\rangle \\
&\ge& \left(1-\frac{1-\lambda}{3(1+\lambda)}\right)\sum_{i\ne j} \frac{\alpha_i^2+\alpha_j^2}{2} \left(k(2j-1, 2i)+k(2j-1, 2i-1)+k(2j, 2i)+k(2j, 2i-1)\right)\\
&-&\frac{2(1-\lambda)}{3(1+\lambda)} (1-\beta)\sum_{i\ne j}|\alpha_i||\alpha_j|\left(k(2i, 2i-1)+k(2j, 2j-1)\right).
\end{eqnarray*}
Also observe that $k(2i, 2i-1)\le \frac{1}{\lambda}(k(2j, 2i)+k(2j-1, 2i))$ and $2k(2j, 2j-1)\le \frac{1}{\lambda}(k(2j-1, 2i-1)+k(2j, 2i-1))$.
The positivity of $\langle (\Ric\wedge \id -\Rm)(\Omega), \Omega\rangle$ holds if
$$
\left(1-\frac{1-\lambda}{3(1+\lambda)}\right)\lambda >\frac{1-\lambda}{3(1+\lambda)}(1- \beta)
$$
which, after plugging in the expression of $\beta$,  is equivalent to
$$
(8m-5)\lambda^2+(6m-3)\lambda -2(m-1)>0.
$$
It is a simple matter to check that if $\lambda>\frac{2(m-1)}{8m-5}$  the above inequality holds up. This completes the proof of the theorem for $\lambda>\frac{2(m-1)}{8m-5}$.

For the case that $\lambda\ge \frac{2(m-1)}{8m-5}(=\frac{n-3}{4n-9})$, observe that the above argument shows that $\Omega$ is parallel (after integration by parts on the manifold). On the other hand, since the manifold has positive sectional curvature at some point and it is of odd dimension its holonomy group has to be $\SO(n)$. The fact that $\Omega$ is parallel implies that as an element in $\mathfrak{so}(n)$ it is fixed by the conjugate action of $\SO(n)$. This is impossible unless $\Omega=0$.\qed

Corollary \ref{berger-5d} is a simple consequence of the vanishing of the second cohomology and the fact that a positively curved compact manifold has finite fundamental group.  We are also aware of the fact that in the case of sectional curvature pinching, Berger did push the constant further down to $4/23$ for $5$-manifolds \cite{B5d}.

\begin{remark} From the proof it is easy to see that the nonnegativity of the complex sectional curvature implies the nonnegativity of Bochner-Weizenb\"ock curvature on two forms, namely $\Ric\wedge \id -\Rm\ge 0$.
\end{remark}

\section{Pinching theorems below one quarter}

First we remark that  Corollary \ref{weak-fqp} carries over to the orbifolds (due to
Theorem \ref{thm:RF-con} in Section 3  and  Proposition 5.2 of \cite{BW}).

\begin{thm} In dimensions above 2  a
 compact Riemannian orbifold with 1/4-pinched flag curvature either admits a metric
of constant sectional curvature or it is isometric to the quotient of a rank 1
symmetric space by a finite group action.
\end{thm}

This generalization is interesting since in even dimensions
 the situation changes dramatically if
one relaxes the assumption:
\begin{prop}\label{prop: exa} In each even dimension $2n>2$ and for each $\lambda<1/4$
there is a compact Riemannian orbifold with sectional curvature $\lambda<K<1$ which
is not given as the quotient of a manifold by a finite group action.
\end{prop}

\noindent {\sl Proof.}
 Consider on $\mathbb{S}^{2n+1}$ (viewed as $(a_1, \ldots, a_{n+1})\in \C^{n+1}$ with $\sum_{i=1}^{n+1} |a_i|^2=1$),  the $\mathsf{S}^1$-action given by
$z(a_1,\ldots,a_{n+1})= (z^{k+1}a_1,z^ka_2,\ldots,z^ka_{n+1})$.
The quotient $M_k\doteqdot\mathbb{S}^{2n+1}/\mathsf{S}^1$ is an orbifold
-- a weighted complex projective space.
If we endow $M_k$ with the induced metric it is straightforward to check
that as $k\to \infty$ the lower and upper curvature bounds of $M_k$ converge
to the lower and upper curvature bounds of $\CP^n$.

On the other hand  $M_k$ is not covered by a manifold. In fact,
consider the frame bundle $F$ corresponding to the
horizontal distribution on $\Sph^{2n+1}$.
The action of $\mathsf{S}^1$ induces a free action on $F$ and the quotient
$\bar{F}=F/\mathsf{S}^1$
can be naturally identified with the frame bundle of the orbifold $M_k$.
Thus the fundamental group of the frame bundle of $M_k$ has at most two elements and
$M_k$ is not a finite quotient of  a manifold.\qed

The above proposition shows that in general the classification of nearly one quarter pinched manifolds of [PT]
can not easily to orbifolds, but see Theorem~\ref{thm: even orbi}.
 In odd dimensions the situation is quite different and we do have the following generalization.

\begin{thm} \label{orbi1}For each constant $C$ and each odd dimension $n$ there is an $\eps>0$
such that the following holds.
Let $(M^n, g)$ be a nonnegatively curved Riemannian orbifold  with $\tfrac{1-\eps}{4}$-pinched flag curvature and  scalar curvature satisfying
$1\le\Scal \le C$. Then $M$ admits a metric of constant curvature.
\end{thm}

\noindent {\sl Proof.}
 First notice that the assumptions  imply the existence
of positive constants $c_1,c_2$ such that the sectional curvature $K$ of $M$ is globally pinched with
$c_1\le K\le c_2$. In order to avoid too much technical difficulties we will
frequently replace the orbifold $M$ by the frame bundle $F$ of the orbifold.
Recall that $F$ is a manifold which can be endowed with a natural connection metric.

The next step  is to show that  the universal cover $\tilde{F}$ of the frame bundle $F$ of the orbifold $M$
has finite second homology.

We may assume that $\tfrac{1-\eps}{4}> \tfrac{n-3}{4n-9}$.
By Theorem \ref{odd-berger} this implies that the Bochner-Weitzenb\"ock operator of $M$ on two forms (namely $\Rm \# \I$) is positive.
We only need to check that the frame bundle $F$  has a metric such that the Bochner-Weitzenb\"ock operator
on two forms is
positive. We endow $F$ with the connection metric  and shrink the fiber $\SO(n)$ by a small factor $\lambda$.
Let $g$ denote the bi-invariant metric on $\SO(n)$.
 It is easy to see that the difference of the curvature tensor
of $(F,g_\lambda)$ and the curvature tensor of the product  $M\times (\SO(n),\lambda^2 g)$
converges to zero.
Since the Bochner-Weitzenb\"ock operators of $M$ and $\SO(n)$ are positive it is easy to see that
the smallest eigenvalue of the Bochner-Weitzenb\"ock operator of
$M\times (\SO(n),\lambda^2 g)$  increases if  $\lambda$ decreases.
This shows that $F$ and its universal cover have finite second homology, by Theorem \ref{odd-berger}.

Next we run the Ricci flow on $M$.
We want to apply a dynamical version of the maximum principle \cite{CL, BW1}.
We assume $\eps<1/2$.
For a large constant $C_2$ and $C_3>>C_2$ and $t\in[0,1/C_3^2]$
we define a set $S(t)\subset S^2_B(\so(n))$ as follows. Let
$S(t)$ be the set of  all algebraic curvature operators
which have $\tfrac{1-\eps-C_3t}{4}$-pinched flag curvature
and whose scalar curvature satisfies $1\le \Scal \le C+C_2t$.
Furthermore we require that for any $\Rm\in S(t)$
the curvature operator $(\Rm+ C_3\eps\cdot e^{C_3^2t}\I)$ has nonnegative
complex sectional curvature.

It is straightforward to check that one can choose $C_2$ and $C_3$ independent of
$\eps$ such that the family $S(t)$ is invariant under the Ricci flow ODE,
\[
\frac{d}{dt} \Rm=2(\Rm^2+\Rm^\#)
\]
that is if $\Rm(t)$ is a solution to the ODE and $\Rm(t')\in S(t')$, then
$\Rm(t)\in S(t)$ for $t\in[t',t_0]$. For example by ODE it is easy to see that there exists $t_1=t_1(C, n)$ such that $|\Rm|\le C'(C)$ on $[0, t_1]$ for some $C'$ depending only on $C$ and dimension $n$. This gives an upper bound $C''(C)$ on the slope $|\frac{d}{dt}\Rm|$. From this it follows easily that $1\le \Scal \le C+C_2t$ is preserved. Similarly one can show that  $\tfrac{1-\eps-C_3t}{4}$-pinched flag curvature is preserved.
To see that the condition on almost nonnegative complex sectional curvature is preserved
one has to use the fact that the nonnegative complex sectional curvature is an invariant condition and
the that the ODE is locally Lipschitz. More precisely, by \cite{BS} (as well as \cite{NW}), if
 $\widehat \Rm= \Rm+ C_3\eps\cdot e^{C_3t}\I$ has nonnegative complex sectional curvature,
then $\widehat{\Rm}^2+\widehat{\Rm}^\#$ is contained in the tangent cone of
complex nonnegative curvature operators.
Since $\frac{d}{dt} \widehat \Rm=\Rm^2+\Rm^{\#}+C^3_3\eps\cdot e^{C_3t}\I$, is greater than $\widehat{\Rm}^2 +\widehat{\Rm}^{\#}$  for large  $C_3$, the same holds for $\frac{d}{dt} \widehat \Rm$.
 By the dynamical version of Hamilton's maximum principle (see, for example, Section 1 of  \cite{BW1})
it follows that the Ricci flow on $M$ exists up to the  time
$t_0$ and $(M,g_t)$ satisfies the curvature condition $S(t)$.

We shall prove more generally that there exists a constant $\eps$
such that the normalized Ricci flow evolves to a constant curvature limit metric.
We argue by contradiction and consider
 a sequence of orbifolds $(M_k,g_k)$ with
$\tfrac{1-\eps_k}{4}$-flag pinching, scalar curvature $1\le \Scal \le C$
and $\eps_k\to 0$ such that for each element in the sequence
the normalized Ricci flow does not converge to a constant curvature limit metric.

Without loss of generality $(M_k,g_k)$ can not be written as a nontrivial
quotient of another orbifold, since otherwise we may replace $M$ by a cover.
This in turn implies that the fundamental group of the frame bundle $F_k$ of $(M,g_k)$
has at most two elements.

We next want to rule out collapsing.
For each fixed $t\in [t_0/4,t_0]$ we can use Shi's estimate to see that
$(M_k,g_k(t))$ has a priori bounds on all derivatives of the curvature tensor.
We now look at the frame bundle $F_k $ of $(M_k,g_k(t))$
with the induced connection metric.
We now rescale the fibers of the frame bundle by a small factor  $\lambda$
independent of $k$, such that there are constants $d_1,d_2>0$ for which
  the  sectional curvature of  this metric on $F_k$  is bounded above by $d_2$
and the Ricci curvature is bounded below by $d_1$.

Since $F_k$ has finite second homology we can employ one of the main theorems of
 \cite{PTu} (see also Theorem 0.2 of \cite{FR})  which asserts that the injectivity radius of $F_k$ is bounded
  from below by  a priori constant depending only on $d_1$ and $d_2$. Now by the
  compactness theorem (cf. \cite{Hamilton-com}, noting that we do not need the compactness result  of  \cite{Lu} on orbifolds),
after passing to a subsequence we obtain  a limit manifold $F$
endowed with a smooth family of metrics $g(t)$ with $t\in [t_0/4,t_0]$ and with
almost free isometric $\SO(n)$-action.
The quotient orbifolds $(M,\bar{g}(t))=(F,g(t))/\SO(n)$ have nonnegative complex sectional
curvature, positive sectional curvature  and $\bar{g}(t)$ is a solution to the Ricci flow.

The strong maximum principle of \cite{BS2} implies that  that either the complex sectional curvature
is positive for $t>t_0/4$ or the orbifold has non-generic holonomy.
Since the sectional curvature is positive and the orbifold is odd dimensional, the holonomy
is generic and it follows that $(M,\bar{g}(t))$ has positive complex sectional curvature for $t>t_0/4$.
This in turn implies that  $(M,g_k(t_0/2))$
has positive complex sectional curvature for infinitely many $k$. This is a contradiction since it implies that
the normalized Ricci flow converges to a constant curvature limit metric for this infinite subsequence.
\qed

\begin{remark}
\begin{enumerate}
\item[a)]The above proof does not give effective bounds on $\eps$.
There is on the other hand a different proof which does give effective bounds.
This alternative proof uses that the Ricci flow ODE behaves somewhat better in odd
dimensions than in even dimensions. In fact with a bit of work one can show
that for an algebraic nonvanishing curvature operator $R\in S^2_B(\so(2n+1))$
with quarter-pinched flag curvature there exists an $\eps>0$ such that
$R+t(R^2+R^\#)$ has positive complex sectional curvature for all $t\in (0,\eps]$.
Knowing this fact one can just apply a dynamical version of the maximum principle
to see that for nearly 1/4-pinched operators, the Ricci flow ODE
pinches towards positive complex sectional curvature.
It would be interesting to see if one can modify this approach such that $\eps$
becomes independent of $C$.
\item[b)]Recall that by a theorem of Abresch and Meyer \cite{AM},
any simply connected odd-dimensional manifold with sectional curvature $K$ satisfying
$\tfrac{1}{4(1+10^{-6})^2}\le K\le 1$ is homeomorphic to a sphere.
An obvious question arises whether or not one can improve the conclusion to diffeomorphism.
\item[c)] For even dimensional manifolds, since the collapsing can not happen, the conclusion of
Petersen and Tao's result still holds under the assumption of Theorem \ref{orbi1}. On the other hand, the argument of the proof above can be sharpen to show that the weighted complex projective spaces are the only exceptions to a similar statement for even dimensional orbifolds. See the theorem below.
\end{enumerate}
\end{remark}

The following theorem shows that the examples constructed in the proof of Proposition~\ref{prop: exa}
are essentially the only additional examples that occur in even dimensions.

\begin{thm} \label{orbi2}\label{thm: even orbi} Given $n$ and $C$ there is an $\eps>0$ such that
any Riemannian orbifold with $\tfrac{1-\eps}{4}$-flag pinching
and scalar curvature satisfying $1\le \Scal \le C$ is either diffeomorphic
to the quotient of rank one symmetric space by a finite isometric group action or
 it is diffeomorphic
to the quotient of a weighted complex projective space by a finite group action.
\end{thm}

The proof requires the following lemma which can be proved
along the same lines as in   the proof of Theorem \ref{main2} (more precisely, the last part of the proof of Theorem \ref{sampson-1}).
\begin{lemma}\label{lem: berger} A quaternionic K\"ahler curvature operator
with quarter pinched flag curvature is a multiple of the curvature operator of $\HP^n$.
\end{lemma}

\begin{proof}[Proof of the theorem.]
We argue  by contradiction.   Consider a sequence of orbifolds $M_k$ satisfying
$\tfrac{1-\eps_k}{4}$-flag pinching,  $1\le \Scal \le C$
with $\eps_k\to 0$, such that each orbifold voliates the conclusion of the theorem.
Similarly to  Theorem \ref{orbi1},  we can find constant $C_2$ and $C_3$ such that the Ricci flow
on $M_k$ exists on $[0, t_0]$ with $t_0=1/C_3^2$ and that $(M_k,g_t)$
satisfies the curvature condition $S(t)$, where $S(t)$ is defined as before.

Recall that geodesics are well defined in the orbifold.
For any regular point $p_k\in M_k$ we consider the exponential map
$\exp\colon T_pM_k\rightarrow M_k$ and pull back the
metric $g_t$ of $M_k$ to the tangent space.
Since we have global curvature bounds, there exists
$r>0$ such that  this metric is nondegenerate
on $B_r(0_{p_k})\subset T_{p_k}M_k$.

We identify this ball with $B_r(0)\subset \R^n$ using a linear
isometry and denote
by  $g_k(p_k,t)$ the induced  Riemannian metric on $B_r(0)$.
By Shi's estimates we have a priori bounds on all the derivatives of the
metric which are independent of $k$ and just depend on a lower bound for
$t\in (0,t_0]$.

We will consider all possible limit metrics  $g(t)$ for all
 convergent subsequences $g_k(p_k,t)$ and all choices $p_k\in M_k$.
We can assume that any limit metric comes with the parameter $t$.
It is then clear that modulo local diffeomorphism $g(t)$ is a solution to the Ricci
flow.
All of the limit metrics have nonnegative complex sectional curvature.
The strong maximum principle can be utilized to see
that the complex sectional curvature is positive unless the limit metric
 has non-generic holonomy. Since the sectional
curvature is positive the only possible non-generic holonomies are
K\"ahler, quaternionic K\"ahler, or $\Spin(9)$-holonomy.
It has been known that the $\Spin(9)$-holonomy implies that
the limit metric is locally isometric to $\CaP^2$.
Below we subdivide the rest of the proof into three cases, which amounts to $\Sph^n$, $\HP^{n/4}$ or $\CaP^2$, and the weighted complex projective spaces.

{\bf Case 1.} There is a sequence of points $p_k\in M_k$ such that
a subsequence  of $g_k(p_k,t)$ converges to a limit metric with generic
holonomy.

After passing to a subsequence we may assume $g_k(p_k,t)$ itself converges.
As explained above the limit metric must have  positive complex sectional curvature.
This in turn implies that we can  find an $\delta>0$ such that the  complex sectional
of $B_{3r/4}(p_k)\subset M_k$ is bounded below by $\delta>0$, for large enough $k$.
Notice that for all $q_k\in B_{5/4r}(p_k)$ the  ball
of radius $r$ around $q_k$ contains a point, on which all complex sectional curvature
$>\delta$. This in turn implies that any convergent subsequence
$g_k(q_k,t)$ has generic holonomy and thus positive complex sectional curvature everywhere.
This in turn shows that there exists a $\delta_2>0$
such that the complex sectional curvature
on the ball $B_{2r}(p_k)$ is bounded below by $\delta'$ for $k$ sufficiently large.
Since we have a priori bounds on the diameter of all orbifolds
a finite iteration of this argument shows that there exists $k_0$ such that
 the complex sectional curvature of $M_k$ is positive for $k\ge k_0$, and thus we get a contradiction.

From now on we assume that any limit metric has  non generic holonomy.

{\bf Claim 1:} There are no two sequences
$p_k, q_k\in M_k$ such that the metrics $g_k(p_k,t)$ and $g_k(q_k,t)$
converge to limit metrics with different holonomy groups.

This follows from Case 1.
In fact, since the diameters of the orbifolds are bounded above
we can choose an fixed integer $l$ and points $p_k^1,\ldots,p_k^l$
satisfying $d(p_k^i,p_k^{i+1})\le r/2$ , $p_k^1=p_k$ and $p_k^l=q_k$.

After passing to subsequence we may assume that $g_k(p_k^j,t)$
converges to a limit metric $g_j$.
Suppose now that the holonomies of the limit metrics $g_j$
and $g_{j+1}$ are different and both not generic.
Suppose for example that $g_j$ is K\"ahler and $g_{j+1}$ is quaternionic
K\"ahler.  Then all the curvature operator  of
$g_k(p_k^j,t)$ are close to being  quaternionic
K\"ahler. Since $B_r(p_k^{j})$ and $B_k(p_k^{j+1})$ have a large intersection
a fixed portion of points in $B_r(0)$ have curvature operators
which are close to ones having quaternionic K\"ahler with respect to the metric
$g_k(p_k^{j},t)$. This in turn shows  for certain points
in $B_r(0)$ the curvature operator (of the limit metric $g_j$) is both K\"ahler and quaternionic K\"ahler.
But this is impossible since such a curvature operator would have a vanishing Ricci tensor, a contradiction as the limit metric has positive sectional curvature.

Hence we below can further assume that all limit metrics have the same holonomy.

{\bf Claim 2:} The limit metric $g(t)$ on $(0, t_0]$
must be locally isometric to a rank one symmetric space.

We only need to consider the case that the limit metric has the holomony of a K\"ahler or quaternionic K\"ahler manifold since $\Spin(9)$-holonomy implies that
the limit metric is locally isometric to $\CaP^2$.
  We only prove it for the case that the limit metric
is K\"ahler,
since in view of Lemma~\ref{lem: berger} the argument for quaternionic K\"ahler case is analogous (easier).
Notice that for $t\to 0$ the limit metrics have
flag pinching constant converging to one quarter which in turn implies that
the curvature operator converges to the curvature operator of $\CP^{n/2}$.
The main idea is to establish a maximum principle for the collection
of all limit metrics.

We consider for each $t$ the supremum $d(t)$ over the following set:
consider for all limit metrics $g(t)$ and all points $x\in B_r(0)$
the distance of the curvature operator $\Rm_{g(t)}(x)$ to nonnegative multiples
of the curvature operator of $\CP^{n/2}$ (resp. of $\HP^{n/4}$).
We claim that the supremum of all these distance is actually attained.

In fact to construct a limit metric where the maximum is attained
we can argue as follows: Choose for each $k$ a point $q_k$
such that the curvature operator $\Rm_{g_k(t)}(q_k)$ has maximal distance $d_k(t)$ to
the O$(n)$-invariant subset of multiples of
the curvature operator of $\CP^{n/2}$.
We choose a regular point $p_k$ with $d(p_k,q_k)<r/2$
and pass to subsequence such that $d_k(t)$ converges to the supremum
and $g_k(p_k,t)$ converges to a limit metric $g(t)$.
Clearly the above supremum is now attained at some point in the closure
of $B_{r/2}(0) \subset B_r(0)$ with this limit metric $g(t)$.

Since the supremum $d(t)$ is attained and $g(t)$ is modulo local diffeomorphisms
a solution to the Ricci flow.
We can now use the maximum principle to derive
that $\lim_{h\to 0_+}\tfrac{d(t)-d(t-h)}{h}\le C_4d(t)$
holds for some universal constant $C_4$.
In fact the space of K\"ahler curvature opertors having distance $\le d(t)$
to muliples of the curvature operators of $\CP^n$ form a convex sets.
Since the Ricci flow ODE leaves the multiples of $\CP^n$
invariant the inequality follows from the dynamical maximum principle
and the fact that $\Rm^2+\Rm^\#$ is locally Lipschitz.

Since $d(t)\to 0$ as $t\to 0$ this inequality implies $d(t)\equiv 0$.

Now we can restrict to the cases that all limits are locally isometric to a rank one symmetric spaces.

{\bf Case 2.} All limit metrics are locally isometric to $\HP^{n/4}$ or $\CaP^2$
up to scaling.

This implies that all limit metrics have a Bochner-Weitzenb\"ock operator
which is positive on two forms. Therefore, the Bochner-Weitzenb\"ock operator
of $M_k$ on two forms is positive for sufficiently large $k$.
As in the odd dimensional case, it follows that
 the universal cover $F_k$ of the frame bundle of $M_k$ has finite
second homology. As before we deduce that $F_k$ has a priori bound
on its injectivity radius and we can assume that $F_k$ converges
to a limit manifold $F$ endowed with a continuous family of metrics.
The quotient $F/\Spin(n)$ is locally isometric to $\HP^{n/4}$  (or $\CaP^2$) and since
$F$ is simply connected it is globally isometric.
Moreover, after passing to a subsequence
 $F_k$ endowed with the action of its isometry group is equivariantly diffeomorphic
to $F$ endowed with a action of a subgroup of its isometry group.

This  shows that $M_k$ is diffeomorphic to a  quotient
of $\HP^{n/4}$ or $\CaP^2$ by a finite isometric group action.

{\bf Case 3.} All limit metrics are locally isometric
to $\CP^{n/2}$ up to scaling.

We replace $M_k$ by a finite cover if needed. It is easy to get a
contradiction if the frame bundle $F_k$ of $M_k$
 has finite second homology, since then noncollapsing
would imply convergence to the frame bundle of $\CP^{n/2}$ which has infinite second homology. Thus we assume that there exists  a nonzero harmonic two form $\theta$ on $F_k$.
It is easy to see that $\theta$ must be equivariant with respect to the $\SO(n)$-action
on $F_k$. In fact otherwise there would exists a non-zero  harmonic
two form $\theta$ on $F_k$ which is perpendicular to any pulled back
2 form on $M_k$. For such a form, it is not hard to check that the $L^2$-norm
of $\nabla \theta$ is bounded below by $\alpha$ times the $L^2$ norm of $\theta$,
where $\alpha>0$ is independent of $k$.
Moreover the number $\alpha$ increases if we scale down the the fibers by a constant.
This is a contradiction  since  the smallest eigenvalue of the Bochner-Weitzenb\"ock operator of $F_k$ converges to $0$ if we scale down the fibers by small factor $\lambda_k\to 0$.
 Hence $\theta$ is the pull back
of some two form $\eta_k$ on $M_k$.

For each $\lambda>0$ we scale down the fibers of the frame bundle $F_k$ by
$\lambda$ and let $\theta_{k,\lambda}$ denote the harmonic two form
with respect to this metric representing a fixed cohomology class.
The above argument shows that $\theta_{k,\lambda}$ is the pull back
of a two form $\eta_{k,\lambda}$ on the orbifold $(M_k,g_k)$.
Clearly $ \eta_{k,\lambda}$ satisfies an ellipic equation and its straightforward to check
that $\eta_{k,\lambda}$ converges to a harmonic two form
$\omega_k$ on $(M_k,g_k)$ for $\lambda\to 0$.

We normalize $\omega_k$ to have $L^2$-norm equal to $\vol(M_k)$.
We may assume  that the pull back $c\cdot pr^* \omega_k$ of $c\omega_k$ to $F_k$ is
a primitive integral class for some $c=c(k)>0$.
We consider the $S^1$ bundle over $F_k$ whose Euler class is $c\omega_k$ and we choose
a connection with curvature $c\cdot \omega_k$.

Since the curvature when restricted to a  $\SO(n)$-fiber vanishes,
the $\SO(n)$-action naturally extends to an isometric $\SO(n)$-action of the total
space of the $S^1$ bundle. By dividing out the $\SO(n)$-action
we obtain an orbifold bundle $S^1\rightarrow S_k\rightarrow M_k$ over $M_k$
with a connection metric whose curvature is given by $c\omega_k$.
The idea is to show that one can scale the $S^1$-fibers such that
the curvature pinching constant of the total space approaches to $1$ as $k\to \infty$.

This  implies that $S_k$ endowed with the $S^1$ action is equivariantly diffeomorphic
to a space form endowed with a linear $S^1$-action and hence the result.

In order to show that the total space has nearly constant curvature
for the right choice of scaling of the fibers it suffices to prove
$\|\nabla \omega_k\|_{L^\infty}\to 0$.
Because this implies that the curvature of the total space $S_k$
approaches the curvature of $S^1\rightarrow S^{2n+1}\rightarrow \CP^n$.

 Integrating the harmonic form gives
\[
\int_{M_k}\int_{B_r(0_p)}\|\nabla \omega_k\|^2_{\exp(v)}=
\vol(B_r(0))\int_{M_k}\|\nabla\omega_k\|^2(p)\le \vol(B_r(p))h_k\vol(M_k)
\]
where $-h_k\to 0 $, as $k\to \infty$,  is the smallest eigenvalue of the Bochner-Weitzenb\"ock operator of $(M_k,g(t))$ on two forms.

For any $q_k\in M_k$ we can find a regular point $p_k$
with $d(p_k,q_k)<r/2$ and
\[
\int_{B_r(0_{p_k})}\|\nabla \omega_k\|^2_{\exp(v)}\le h_k\vol(B_r(0))\tfrac{\vol(M_k)}{\vol(B_{r/2}(q_k))}\doteqdot j_k
\]
with $j_k\to 0$, where we have used that the ratio $\tfrac{\vol(M_k)}{\vol(B_{r/2}(q_k))}$ is bounded.

%

Since $\exp^*\omega_k$ satisfies an elliptic equation
with respect to a metric for which we have a priori bounds  on all derivatives
we deduce that $\exp^* \omega_k$ converges to a form which is parallel
with respect to the limit metric $g(t)$.
This in turn shows that on the orbifold  $\|\nabla \omega_k\|_{L^\infty}\to 0$.

\end{proof}

\section*{Acknowledgments.} { The first author thanks Nolan Wallach for some helpful discussions.  His research is partially supported by a NSF grant DMS-0805834. The paper was written  when the second author visited UCSD during February--March of 2009. The second author thanks the mathematical department of UCSD for its hospitality. }

\bibliographystyle{amsalpha}

\begin{thebibliography}{A}

\bibitem[AM]{AM} U. Abresch and W.-T. Meyer, \textit{ A sphere theorem with a pinching constant below ${1\over4}$.}  J. Differential Geom.  \textbf{44}(1996),  no. 2, 214--261.

\bibitem[AN]{AN} B. Andrews and H. Nguyen, \textit{ Four-manifolds with $1/4$-pinched flag curvature.} Preprint.


\bibitem[B1]{Berger} M. Berger, \textit{ Sur quelques vari\'et\'es riemanniennes suffisamment pinc\'ees.} Bull. Soc. Math. France \textbf{88} (1960), 57--71.

\bibitem[B2]{B5d} M, Berger, \textit{ Sur les varites $(4/23)$-pinces de dimension 5.} C. R. Acad. Sci. Paris, \textbf{257}(1963), 4122--4125.


\bibitem[BW1]{BW1} C. B\"ohm and B. Wilking, \textit{ Nonnegatively curved manifolds with finite fundamental groups admits metrics with positive Ricci curvature.}   Geom. Funct. Anal. \textbf{17}(2007), 665--681.


\bibitem[BW2]{BW} C. B\"ohm and B.  Wilking,  \textit{ Manifolds with positive curvature operators are space forms.}  Ann. of Math. (2)  \textbf{167}(2008),  no. 3, 1079--1097.

\bibitem[BS1]{BS} S. Brendle and R. Schoen, \textit{ Manifolds with
$1/4$-pinched curvature are space forms.} Jour. Amer. Math. Soc.  \textbf{22}(2009),  no. 1, 287--307.


\bibitem[BS2]{BS2}S. Brendle and R. Schoen, \textit{  Classification of manifolds with weakly $1/4$-pinched curvatures.}  Acta Math.  \textbf{200}(2008),  no. 1, 1--13.

\bibitem[Ch]{Ch} H. Chen, \textit{ Pointwise $1/4$-pinching 4 manifolds.} Ann. Global Geom. \textbf{9}(1991), 161--176.

\bibitem[CL]{CL} B. Chow and P. Lu, \textit{ The maximum principle for systems of parabolic equations subject to an avoidance set.}  Pacific J. Math.  \textbf{214}(2004),  no. 2, 201--222.

\bibitem[CY]{CY} B.
Chow and D. Yang, \textit{ Rigidity of nonnegatively curved compact quaternionic-K\"ahler manifolds.}  J. Differential Geom.  \textbf{29}(1989),  no. 2, 361--372.

\bibitem[FR]{FR} F. Fang and X. Rong, \textit{ Positive pinching, volume and second Betti number.}  Geom. Funct. Anal.  \textbf{9}(1999),  no. 4, 641--674.


\bibitem[Gr]{Gr} M.
Gromov, \textit{ Positive curvature, macroscopic dimension, spectral gaps and higher signatures.}  Functional analysis on the eve of the 21st century, Vol. II (New Brunswick, NJ, 1993),  1--213, Progr. Math., 132, Birkhäuser Boston, Boston, MA, 1996.




\bibitem[H1]{H82} R. Hamilton, \textit{ Three-manifolds with positive Ricci curvature.} J. Differential Geom. \textbf{17}(1982), 255--306.

\bibitem[H2]{H86} R. Hamilton, \textit{ Four-manifolds with positive curvature operator.}  J. Differenital.
Geom. \textbf{24}(1986 ), 153--179.

\bibitem[H3]{Hamilton-com} R. Hamilton, \textit{
A compactness property for solutions of the Ricci flow.}  Amer. J. Math.  \textbf{117}(1995),  no. 3, 545--572.

\bibitem[He]{He} L.
Hern\'andez, \textit{ K\"ahler manifolds and $1/4$-pinching.}  Duke Math. J.  \textbf{62}(1991),  no. 3, 601--611.

\bibitem[L]{Lu} P. Lu, \textit{A compactness property for solutions of the Ricci flow on orbifolds.}  Amer. J. Math.  \textbf{123}(2001),  no. 6, 1103--1134.

\bibitem[Ng]{Ng} H. Nguyen, \textit{ Invariant curvature cones and the Ricci flow.} Ph. D thesis. Australian National University, 2007.


\bibitem[NW]{NW} L. Ni and J. Wolfson, \textit{ Positive complex sectional curvature, Ricci flow and the differential shpere theorem. } Unpublished manuscript.




\bibitem[PT]{peter-tao}  P. Petersen and T. Tao, \textit{ Classification of almost quarter-pinched manifolds.} Preprint.

\bibitem[PTu]{PTu} A. Petrunin and W. Tuschmann, \textit{ Diffeomorphism finiteness, positive pinching, and second homology.}  Geom. Funct. Anal. \textbf{9}(1999), 736--774.

\bibitem[Sa]{Sa} J. H.
Sampson, \textit{ Harmonic maps in K\"ahler geometry.}  Harmonic mappings and minimal immersions (Montecatini, 1984),  193--205, Lecture Notes in Math., 1161, Springer, Berlin, 1985.

\bibitem[Si]{Siu}Y.-T. Siu, \textit{  The complex-analyticity of harmonic maps and the strong rigidity of compact K\"ahler manifolds.}  Ann. of Math. (2)  \textbf{112}(1980), no. 1, 73--111.

\bibitem[W]{W} B. Wilking, \textit{ Nonnegatively and positively curved manifolds.}  Surveys in differential geometry. Vol. XI,  25--62, Surv. Differ. Geom., 11, Int. Press, Somerville, MA, 2007.


\bibitem[YZ]{YZ} S.T. Yau and F. Zheng, \textit{Negatively
$1/4$-pinched Riemannian metric on a compact K\"ahler manifold.}
Invent. Math. \textbf{103}(1991), 527--535.

\end{thebibliography}

{\sc  Addresses:}

{\sc  Lei Ni},
 Department of Mathematics, University of California at San Diego, La Jolla, CA 92093, USA


email: lni@math.ucsd.edu

{\sc Burkhard Wilking},
Mathematisches Institut, University of M\"unster,  Einsteinstrasse 62, 48149 M\"unster, Germany

email: wilking@math.uni-muenster.de

\end{document}